\documentclass[12pt]{article}
\usepackage{amssymb,amsmath,latexsym}

\DeclareMathOperator{\Mod}{mod}

\def\B{{\cal B}}
\def\FC{{\mathfrak C}}
\def\CC{{\cal C}}
\def\cycp{\vec{C}^p}
\def\cycn{\vec{C}^n}

\def\E{{\cal E}}
\def\F{{\cal F}}
\def\R{{\cal R}}
\def\S{{\cal S}}
\def\lg{\langle}
\def\rg{\rangle}

\def\proof{{\bf Proof}.\ }
\def\myskip{\vspace{2mm}}
\def\bull{\hfill\vrule height .9ex width .8ex depth -.1ex\myskip}

\newtheorem{formula}{}[section]

\newtheorem{definition}[formula]{Definition}
\newtheorem{corollary}[formula]{Corollary}

\newtheorem{lemma}[formula]{Lemma}
\newtheorem{theorem}[formula]{Theorem}

\def\thrm{\begin{theorem}}
\def\thrml#1{\begin{theorem}\label{#1}}
\def\ethrm{\end{theorem}}
\def\dfntn{\begin{definition}}
\def\dfntnl#1{\begin{definition}\label{#1}}
\def\edfntn{\end{definition}}
\def\nmrt{\begin{enumerate}}
\def\enmrt{\end{enumerate}}
\def\tm#1{\item[{\rm (#1)}]}
\def\qtn{\begin{equation}}
\def\qtnl#1{\begin{equation}\label{#1}}
\def\eqtn{\end{equation}}
\def\lmm{\begin{lemma}}
\def\lmml#1{\begin{lemma}\label{#1}}
\def\elmm{\end{lemma}}
\def\crllr{\begin{corollary}}
\def\crllrl#1{\begin{corollary}\label{#1}}
\def\ecrllr{\end{corollary}}
\def\css{\begin{cases}}
\def\ecss{\end{cases}}

\title{ \bf{The basis digraphs of $p$-schemes}}
\author{
Ilia Ponomarenko
\thanks{Partially supported by RFFI grants 05-01-00899, NSH-4329.2006.1}\\[-1pt]
\small Corresponding author\\[-3pt]
\small Petersburg Department of V.A.Steklov\\[-3pt]
\small Institute of Mathematics\\[-3pt]
\small Fontanka 27, St. Petersburg 191023, Russia\\[-3pt]
{\tt \small inp@pdmi.ras.ru}\\[-3pt]
\small http://www.pdmi.ras.ru/\~{}inp
\and A. Rahnamai Barghi
\thanks{Partially supported by IASBS, Zanjan- Iran. The author was visiting the
Euler Institute of Mathematics, St. Petersburg, Russia during the time a part
of this paper was written and he thanks the Euler Institute for its hospitality.}\\
\small Institute for Advanced Studies in Basic Sciences (IASBS),\\
\small P.O.Box: 45195-1159, Zanjan,
Iran\\
{\tt \small rahnama@iasbs.ac.ir}\\[-3pt]
\small http://www.iasbs.ac.ir/faculty/rahnama/ }

\date{September 28, 2007}

\begin{document}

\maketitle

\begin{abstract}
It is proved that association schemes with  bipartite basis graphs are
exactly 2-schemes. This result follows from a characterization of $p$-schemes
for an arbitrary prime $p$ in terms of basis digraphs.\myskip

\end{abstract}
\newpage

\section{Introduction}
Nowadays, the theory of association schemes is usually considered as a generalization of the
theory of finite groups. In this sense, $p$-schemes introduced in~\cite{Zi1} for a prime
number~$p$ give a natural analog of $p$-groups that enables us, for example, to get the Sylow
theorem for association schemes~\cite{HMZ}. It is a routine task to extend this notion to
coherent configurations (for short, schemes) the special case of which are association schemes
(for the exact definitions and notations concerning schemes see section~\ref{020507a}). This was
done in~\cite{PB} where algebraic properties of $p$-schemes were studied. In this paper
we focus on combinatorial features of $p$-schemes.\myskip

From a combinatorial point of view an association scheme $\CC$ can be thought as a special partition of a
complete digraph into spanning subdigraphs satisfying certain regularity conditions.
These subdigraphs are called the {\it basis digraphs} of $\CC$; exactly one of them consists in all
loops and is called a {\it reflexive} one. In order to state our main result, we need the
following graph-theoretical notion: given an integer $p>1$ a digraph $\Gamma$ is said to be
{\it cyclically $p$-partite} if its vertex set can be partitioned into $p$ nonempty mutually
disjoint sets in such a way that if a pair $(u,v)$ is an arc of $\Gamma$, where $u$ (resp. $v$)
belongs to $i$-th (resp. $j$-th) set, then $j-i$ is equal to $1$ modulo $p$, see~\cite[p.82]{CDS80}.

\thrml{261106a}
Let $p$ be a prime and $\CC$ an association scheme. Then $\CC$ is a $p$-scheme if and only if
each non-reflexive basis digraph of~$\CC$ is cyclically $p$-partite.
\ethrm

From \cite[Theorem~3.4]{PB} it follows that a characterization of arbitrary $p$-schemes can be
reduced to association scheme case in which Theorem~\ref{261106a} works. Besides, a digraph~$\Gamma$
is $2$-partite if and only if the corresponding undirected loopless graph $\Gamma'$ is
bipartite. (When $\Gamma$ is a basis digraph of a scheme~$\CC$, we say that $\Gamma'$ is the
{\it basis graph} of~$\CC$.) Thus we come to the following characterization of $2$-schemes.

\crllrl{291106}
Let $\CC$ be a scheme. Then $\CC$ is a $2$-scheme if and only if each basis graph of~$\CC$ is
bipartite.
\ecrllr

The proofs of Theorem \ref{261106a} and Corollary \ref{291106} will be given in Section \ref{020507b}.
Section~\ref{020507a} contains notations and definitions concerning schemes. In
Section~\ref{020507c} we prove several results on $p$-schemes which will be used later in the
proof of Theorem \ref{261106a}.
\myskip

{\bf Notations.}
Throughout the paper $V$ denotes a finite set.

By a relation on $V$ we mean any set $R\subseteq V\times V$. The smallest set $X\subseteq V$
such that $R\subseteq X\times X$ is called the support of $R$ and is denoted by $V_R$.
Set $\Delta(V)=\{(v,v): v\in V\}$ to be the diagonal relation on~$V$.

Given $R,S\subseteq V\times V$ we set $RS=\{(u,v)\in V\times V: (u,w)\in R$ and $(w,v)\in S$ for some
$w\in V\}$ and call it the product of~$R$ and~$S$.

Given sets $X,Y\subseteq V$ and a set $\R$ of relations on $V$ we denote by $\R_{X,Y}$
the set of all nonempty relations $R_{X,Y}=R\cap(X\times Y)$ with $R\in\R$. We
write $\R_X$ and $R_X$ instead of $\R_{X,X}$ and $R_{X,X}$ respectively.

By an equivalence $E$ on $V$ we mean an ordinary equivalence relation on a subset of~$V$.
The set of its classes is denoted by~$V/E$. Given $X\subseteq V$ we set $X/E=X/E_X$.
The set of all equivalences on~$V$ is denoted by~$\E_V$.

Given an equivalence $E\in\E_V$ and a set $\R$ of relations on $V$ we denote by $\R_{V/E}$
the set of all nonempty relations $R_{V/E}=\{(X,Y)\in V/E\times V/E:\ R_{X,Y}\neq\emptyset\}$
where $R\in\R$.

By a digraph we mean a pair $\Gamma=(V,R)$ where $R$ is a relation on~$V$.
The digraph is called reflexive if $\Delta(V)\subseteq R$ .

A cycle of length $n$ is the digraph $\cycn=(V,R)$ where $V=\{0,\ldots,n-1\}$ and
$R$ consists of arcs $(i,i+1)$, $i\in V$, with the addition taken modulo~$n$.

\section{Schemes}\label{020507a}
Let $V$ be a finite set and $\R$ a partition of $V\times V$ closed with respect to the permutation
of coordinates. Denote by $\R^*$ the set of all unions of the elements of $\R$. A pair
$\CC=(V,\R)$ is called a {\it coherent configuration}~\cite{H70} or a {\it scheme} on $V$ if
the set $\R^*$ contains the diagonal relation~$\Delta(V)$, and given $R,S,T\in\R$,
the number
$$
c_{R,S}^T=|\{v\in V:\,(u,v)\in R,\ (v,w)\in S\}|
$$
does not depend on the choice of $(u,w)\in T$. The elements of $V$ and $\R$ are called the
{\it points} and the {\it basis relations} of $\CC$ respectively. Two schemes $\CC$ and $\CC'$
are called {\it isomorphic}, $\CC\cong\CC'$, if there exists a bijection between their point
sets which preserves the basis relations.\myskip

A set $X\subseteq V$ is called a {\it fiber} of~$\CC$ if the diagonal relation
$\Delta(X)$ is a basis one. Denote by~$\F$ the set of all fibers. Then
$$
V=\bigcup_{X\in\F}X,\qquad \R=\bigcup_{X,Y\in\F}\R_{X,Y}
$$
where the both unions are disjoint. The scheme $\CC$ is called  {\it homogeneous} or {\it association scheme}
if $|\F|=1$. In this case $\Delta=\Delta(V)$ is a basis relation of it and
\qtnl{030907a}
c_{R,R^T}^\Delta=c_{R^T,R}^\Delta=|R(u)|,\quad R\in\R,\ u\in V,
\eqtn
where $R^T=\{(u,v):\ (v,u)\in R\}$ and $R(u)=\{v\in V: (u,v)\in R\}$. In particular,
the cardinality of the latter set does not depend on~$u$. We denote it by $d(R)$.
Clearly, $|R|=d(R)|V|$ for all $R\in\R$.\myskip

By an {\it equivalence} of the scheme~$\CC$ we mean any element of the set $\E=\E(\CC)=\R^*\cap\E_V$.
Given $E\in\E$ one can construct schemes
$$
\CC_{V/E}=(V/E,\R_{V/E}),\qquad \CC_X=(X,\R_X)
$$
where $X\in V/E$. If the set~$X$ is a fiber of~$\CC$, then obviously $X\times X\in\E$, and hence
the scheme $\CC_X$ is homogeneous.
The equivalence $E\neq\Delta$ is {\it minimal} if no other equivalence in $\E\setminus \{\Delta\}$
is contained in~$E$, and $E\neq V\times V$ is called {\it maximal} if no other equivalence in
$\E\setminus\{V\times V\}$ contains~$E$. The set of all maximal (resp. minimal) equivalences
of~$\CC$ is denoted by~$\E_{max}$ (resp.~$\E_{min}$). A homogeneous scheme $\CC$ on at least
two points is called {\it primitive} if $\E=\{\Delta,V\times V\}$.
Clearly, $E\in\E_{max}$ (resp. $E\in\E_{min}$) if and only if the scheme $\CC_{V/E}$
(resp. $\CC_X$ for some $X\in V/E$) is primitive.\myskip

For a homogeneous scheme $\CC$ the set
\qtnl{300807a}
G=\{R\in\R:\ d(R)=1\}
\eqtn
forms a group with
respect to the product of relations. The identity of this group coincides with $\Delta$.
The order of an element $R\in G$ equals the sum of all numbers $d(S)$ where $S$ is a basis
relation of~$\CC$ contained in the set
\qtnl{030907b}
\langle R \rangle=\bigcup_{i\ge 0}R^i.
\eqtn
A set $\S\subseteq\R$ is a subgroup of $G$ if and only if the
union of all $R\in\S$ belongs to~$\E$. The scheme $\CC$  is called {\it regular}, if
$G=\R$.\myskip

Given $R\in\R^*$ the digraph $\Gamma(\CC,R)=(V_R,R)$ is called the {\it basis digraph} (resp.
the {\it basis graph}) of a scheme~$\CC$, if $R\in\R$ (resp. $R=(S\cup S^T)\setminus\Delta$ for
some $S\in\R$). In particular, any basis graph of $\CC$ is an undirected loopless graph.
From~\cite[p.55]{W76}, it follows that the basis digraph of a homogeneous scheme is strongly
connected if and only if the corresponding basis graph is connected. This implies that the
relation~(\ref{030907b}) is an equivalence of the scheme~$\CC$. It is easy to see that
it is the smallest equivalence on~$V$ containing~$R$.

\section{$p$-schemes}\label{020507c}

Throughout this section $p$ denotes a prime number. A scheme $\CC=(V,\R)$ is called a {\it $p$-scheme}
if the cardinality of any relation $R\in\R$ is a power of $p$ (for more details see~\cite{PB}
a~\cite{Zi1}). The class of all $p$-schemes is denoted by~$\FC_p$.

\thrml{060507c}
Let $\CC\in\FC_p$ be a primitive scheme. Then $\CC$ is a regular and $|V|=p$. In particular,
any non-reflexive basis digraph of~$\CC$ is isomorphic to~$\cycp$.
\ethrm
\proof By the assumption $\CC$ is a homogeneous scheme. Due to (\ref{030907a}) this implies that
$$
|\Delta|=|V|=\sum_{R\in\R}d(R).
$$
Since $\CC\in\FC_p$, both the left-hand side and each summand in the right-hand side are
powers of~$p$. Taking into account that $d(\Delta)=1$ and $|V|\ge 2$ (because of the
primitivity), we conclude that there exists a non-diagonal relation $R\in\R$ such that $d(R)=1$.
By \cite[p.71]{W76} any primitive scheme having such a basis relation $R$ is a regular
scheme on $p$ points.\bull

A special case of the following statement was proved in \cite{PB}.

\thrml{130207}
Let $\CC$ be a homogeneous scheme, $E\in\E$ and $X\in V/E$. Then $\CC\in\FC_p$ if and
only if $\CC_{V/E}\in \FC_p$ and $\CC_X\in\FC_p$.
\ethrm
\proof The necessity follows from the obvious equality
\qtnl{291106b}
|R|=|R_{V/E}|\cdot|R_{X,Y}|,\qquad R\in\R,
\eqtn
where $X,Y\in V/E$ with $R_{X,Y}\ne\emptyset$. Let us prove the sufficiency. Without
loss of generality we may assume that $|V|>1$. Suppose that $E\not\in\E_{min}$. Then there
exists an equivalence $F\in\E_{min}$ such that $F\subsetneq E$. The scheme $\CC'=\CC_{V/F}$
is a homogeneous one and by \cite[Theorem~1.7.6]{Zi1} we have
$$
\CC'_{V'/E'}\cong\CC_{V/E}^{},\qquad \CC'_{X'}\cong(\CC_X^{})_{X/F}^{},
$$
where $V'=V/F$, $E'=E_{V/F}$ and $X'$ is the class of $E'$ such that $X/F=X'$. From the
first part of the proof (for $\CC=\CC_X$ and $E=F_X$) it follows that $(\CC_X)_{X/F}\in\FC_p$. Since
$\CC_{V/E}\in\FC_p$ and $|V'|<|V|$, we conclude by induction that $\CC_{V/F}=\CC'\in\FC_p$.
Thus we can replace $E$ by the equivalence~$F\in\E_{min}$. In this case the scheme $\CC_X$ is a
primitive $p$-scheme. By Theorem~\ref{060507c} it is a regular scheme on~$p$ points. Thus
$\CC\in\FC_p$ by~\cite[Theorem~3.2]{PB}.\bull

A set $X\subseteq V$ is called a {\it block} of a scheme $\CC$ if there exists an equivalence
$E\in\E$ such that $X\in V/E$. Clearly, $V$ is a block; it is called a {\it trivial} one.
Denote by $\B$ the set of all nontrivial blocks of~$\CC$.

\thrml{151106a}
Let $\CC$ be a homogeneous scheme satisfying the following conditions:
\nmrt
\tm{1} $|\E_{max}|\ge 2$,
\tm{2} $\CC_X\in\FC_p$ for all $X\in\B$.
\enmrt
Then $\CC\in\FC_p$.
\ethrm
\proof First suppose that the scheme $\CC$ is regular. Then it suffices to verify that the
group~$G$ defined by formula~(\ref{300807a}) is a $p$-group. However, from condition~(2) it
follows that any proper subgroup of~$G$ is a $p$-group. This implies that $G$ is a $p$-group
unless it is of prime order other than~$p$. Since the latter contradicts to condition~(1), we
are done.\myskip

Suppose that $\CC$ is not regular. By condition (1) there are distinct equivalences $E_1,E_2\in\E_{max}$.
Without loss of generality we can assume that there exists an equivalence $E\in\E_{min}$ such that
\qtnl{291106a}
E\subseteq E_1\cap E_2.
\eqtn
Indeed, if $E_1\cap E_2\ne\Delta$, then one can take as $E$ a minimal equivalence of~$\CC$ contained
in~$E_1\cap E_2$. Otherwise,
\qtnl{060507h}
F_i\cap E_j=\Delta,\qquad \{i,j\}=\{1,2\}.
\eqtn
where $F_i$ is a minimal equivalence of $\CC$ contained in $E_i$, $i=1,2$. From
Theorem~\ref{060507c} (applied for $\CC=\CC_X$ with $X\in V/F_i$) it follows that
$F_1,F_2\subseteq G$. Since $G$ is closed with respect to products of relations, this implies
that it contains the subgroup $F=\lg F_1,F_2\rg$. Moreover, $F\ne V\times V$, for otherwise the
scheme $\CC$ is regular which contradicts to the assumption. So there exists an equivalence
$F'\in\E_{max}$ such that $F\subseteq F'$. We observe that $E_1\ne F'$, for otherwise
$$
F_2\subseteq \lg F_1,F_2\rg=F\subseteq F'=E_1
$$
which contradicts to (\ref{060507h}). Since $F_1\subseteq E_1$ and $F_1\subseteq F\subseteq F'$, this shows that
inclusion~(\ref{291106a}) holds for $E=F_1$ and $E_2=F'$.\myskip

From (\ref{291106a}) it follows that $(E_1)_{V/E}$ and $(E_2)_{V/E}$ are distinct maximal
equivalences of the scheme $\CC'=\CC_{V/E}$. In particular, $|\E'_{max}|\ge 2$ where
$\E'=\E(\CC')$. Besides, any block of $\CC'$ is of the form $X'=X/E$ for some block~$X$ of~$\CC$.
By condition~(2) and Theorem~\ref{130207} this shows that
$$
\CC'_{X'}=(\CC_{V/E})_{X'}\cong(\CC_X)_{X/E}\in\FC_p.
$$
Thus the scheme $\CC'$ satisfies conditions~(1) and~(2). Since, obviously, $|V/E|<|V|$, it
follows by induction that $\CC'\in\FC_p$. By Theorem~\ref{130207} this implies that
$\CC\in\FC_p$ and we are done.\bull

It should be remarked that condition (1) in Theorem~\ref{151106a} is essential. Indeed,
let $\CC$ be the wreath product of a regular scheme on $p$ points by a regular scheme on
$q$ points where $p$ and $q$ are different primes \cite[p.45]{W76}. Then the set $\E_{max}=\E_{min}$
consists of a unique equivalence $E$ such that $\CC_X$ is a regular scheme on $p$ points
for all $X\in V/E$. Thus $\CC$ satisfies condition (2), does not satisfy condition (1) and
is not a $p$-scheme.

\section{Proofs of Theorem \ref{261106a} and Corollary \ref{291106}}\label{020507b}

Throughout this section we fix $n=|V|$ and an integer $p>1$. By the definition given in the
introduction a digraph $\Gamma=(V,R)$ is cyclically $p$-partite if and only if the set $V$ is a
disjoint union of nonempty sets $V_0,\ldots,V_{p-1}$ such that
\qtnl{030507f}
R=\bigcup_{r=0}^{p-1}R_{V_r,V_{r+1}}
\eqtn
with addition taken modulo $p$. In particular, $n\ge p$, and $n=p$ if and only if $\Gamma$ is
isomorphic to a subdigraph of the directed cycle $\cycp$. The following statement shows that
the class of cyclically $p$-partite digraphs is closed with respect to taking a {\it disjoint
union} where under the disjoint union of digraphs $(V_i,R_i)$, $i\in I$, we mean
the digraph $(V,R)$ with $V$ and $R$ being the disjoint unions of $V_i$'s and $R_i$'s
respectively.

\lmml{030507a}
Let $\Gamma$ be a disjoint union of strongly connected digraphs $\Gamma_i$, $i\in I$. Then
$\Gamma$ is cyclically $p$-partite if and only if so is $\Gamma_i$ for all $i$.
\elmm
\proof The sufficiency is clear. To prove the necessity suppose that $\Gamma=(V,R)$ and
$V$ is a disjoint union of nonempty sets $V_0,\ldots,V_{p-1}$ for which equality~(\ref{030507f})
holds. Let us verify that given $i\in I$ the digraph $\Gamma_i=(X,S)$ is cyclically $p$-partite.
We observe that from (\ref{030507f}) it follows that
$S(x)\subseteq R(x)\subseteq V_{r+1}$ for all $x\in V_r\cap X$ and all $r$. On the other
hand, since $\Gamma_i$ is a strongly connected digraph, we also have
$S(x)\ne\emptyset$ for all $x\in X$. Thus
$$
V_r\cap X\ne\emptyset\ \Rightarrow\ V_{r+1}\cap X\ne\emptyset,\qquad r=0,\ldots,p-1,
$$
whence it follows that $V_r\cap X\ne\emptyset$ for all $r$. Since obviously equality
(\ref{030507f}) holds for $R=S$ and $V_r=V_r\cap X$, we conclude that $\Gamma_i$ is cyclically
$p$-partite.\bull

Let $\Gamma=\Gamma(\CC,R)$ be a basis digraph of a homogeneous scheme $\CC$. Then $\Gamma$ is
strongly connected if and only if the graph $(V,R\cup R^T)$ is connected (see Section~\ref{020507a}),
or equivalently $\lg R\rg=V\times V$. This implies that in any case the digraph $\Gamma$ is
disjoint union of digraphs $\Gamma(C_X,R_X)$ where $X$ runs over the classes of the
equivalence $\lg R\rg$. By Lemma~\ref{030507a} this proves the following statement.

\crllrl{060507a}
Let $\CC$ be a homogeneous scheme. Then given a non-diagonal basis relation $R\in\R$ the
digraph $\Gamma(\CC,R)$ is cyclically $p$-partite if and only if so is the digraph
$\Gamma(\CC_X,R_X)$ for all $X\in V/\lg R\rg$.\bull
\ecrllr

{\bf Proof of Theorem~\ref{261106a}.} Let $\CC=(V,\R)$ be a homogeneous scheme. Without
loss of generality we may assume that $n>1$.\myskip

To prove necessity, suppose that $\CC\in\FC_p$, $\Gamma(\CC,R)$ is a non-reflexive
basis digraph of~$\CC$ and $E=\lg R\rg$. If $E\ne V\times V$, then $|X|<n$, and $\CC_X\in\FC_p$
for all $X\in V/E$ (statement~(1) of Theorem~\ref{130207}). By induction this implies that
the digraph $\Gamma(\CC,R_X)$ is cyclically $p$-partite for all $X$, and we are done
by Corollary~\ref{060507a}. Let now $E=V\times V$. Take $F\in\E_{max}$. Then $\CC_{V/F}$ is
a primitive $p$-scheme (Theorem~\ref{130207}) and $R_{V/F}$ is a non-diagonal basis relation
of it. By Theorem~\ref{060507c} this implies that
$$
\Gamma(\CC_{V/F},R_{V/F})\cong\cycp.
$$
Therefore, the equivalence $F$ has $p$ classes, say $V_0,\ldots,V_{p-1}$, and equality~(\ref{030507f})
holds for a suitable numbering of $V_i$'s. Thus the graph $\Gamma(\CC,R)$ is cyclically $p$-partite.
\myskip

To prove the sufficiency, suppose that each non-reflexive basis digraph of $\CC$ is cyclically
$p$-partite.  Then by Corollary~\ref{060507a} so is each non-reflexive basis digraph of $\CC_X$
for all $X\in\B$. By induction this implies that
\qtnl{040507h}
\CC_X\in\FC_p,\qquad X\in\B.
\eqtn
So if $|\E_{max}|\ge 2$, then $\CC\in\FC_p$ by Theorem~\ref{151106a} and we are done. Otherwise,
$\E_{max}=\{F\}$ for some equivalence $F\in\E$. Take a relation $R\in\R$ such that
$R\cap F=\emptyset$. Then $\lg R\rg\not\subseteq F$ and hence $\lg R\rg=V\times V$.
In particular, the digraph $\Gamma=\Gamma(\CC,R)$ is strongly connected. Since it is also
cyclically $p$-partite, there exists an equivalence $E\in\E_V$ with $p$ classes
$V_0,\ldots,V_{p-1}$ for which equality (\ref{030507f}) holds. The strong connectivity
of~$\Gamma$ implies that
$$
(u,v)\in E\ \Leftrightarrow\ d(u,v)\equiv 0\,(\Mod p)
$$
where $d(u,v)$ denotes the distance between $u$ and $v$ in the graph $\Gamma$. It follows
that $E$ is a union of relations $R^{ip}$ where $i$ is a nonnegative integer. Therefore
$E\in\E$. This enables us to define the scheme~$\CC_{V/E}$. Due to~(\ref{030507f}) we have
$$
\Gamma(\CC_{V/E},R_{V/E})\cong\cycp.
$$
So $\CC_{V/E}$ is a scheme on $p$ points having basis relation $R_{V/E}$ with $d(R_{V/E})=1$.
It follows that $\CC_{V/E}\in\FC_p$. Together with~(\ref{040507h}) this shows that
$\CC$ satisfies the sufficiency condition of Theorem~\ref{130207}. Thus $\CC\in\FC_p$.\bull
\myskip

{\bf Proof of Corollary \ref{291106}}. Let $R$ be a basis relation of the scheme~$\CC$. It is
easy to see that a non-reflexive digraph $(V,R)$ is cyclically $2$-partite if and only if the
graph $(V,R\cup R^T)$ is bipartite. Besides, by \cite[Theorem~3.4]{PB}
a scheme $\CC$ is a $p$-scheme if and only if so is the scheme $\CC_X$ for all $X\in\F$. Thus
Theorem~\ref{261106a} implies that $\CC$ is a $2$-scheme if and only if the graph $\Gamma(\CC_X,R\cup R^T)$ is
bipartite for all $X\in\F$ and all $R\in\R_{X,X}\setminus\{\Delta(X)\}$. This proves the sufficiency.
The necessity follows from the fact that given $R\in\R_{X,Y}$ with distinct $X,Y\in\F$,
the graph $\Gamma(\CC,R\cup R^T)$ is bipartite.\hfill \bull

\end{document}